\documentclass[11pt,a4paper,reqno,centertags]{article}
\usepackage{xypic,a4wide,amsmath,amsthm,amssymb,tabularx,calc}

\newtheorem{thm}{Theorem}[section]
\newtheorem{corol}[thm]{Corollary}

\newtheorem{prop}[thm]{Proposition}
\newtheorem{defin}[thm]{Definition}

\theoremstyle{remark}
\newtheorem{rem}[thm]{Remark}
\newtheorem{ex}[thm]{Example}

\newcommand{\qee}{\mbox{\hspace{0.2mm}}\hfill$\triangle$}
\def\prf{\noindent{\textsc{Proof}}\rm\ }
\def\endprf{\ \hfill $\Box$\medskip}

\def\cal{\mathcal}

\def\bC{{\Bbb C}}

\def\bR{{\Bbb R}}
\def\bZ{{\Bbb Z}}
\def\bN{{\Bbb N}}

\def\ddz{\frac{dz}{2\pi i z}}
\def\SC{\mathrm{SC}}
\def\Id{\mathrm{I}}

{\nonumber}

\begin{document}
\begin{center}
 {\LARGE\bf Matrix biorthogonal polynomials on the unit circle and non-Abelian Ablowitz-Ladik hierarchy.}\\[15pt]
 {\sc Mattia Cafasso}\\
 {Universit\'e Catholique de Louvain}\\
 {D\'epartement de Math\'ematiques\\
Chemin du Cyclotron, 2\\
1348 Louvain-la-Neuve\\
Belgium }
\end{center}

\vspace{1cm}

\begin{abstract}
	Adler and van Moerbeke \cite{AVM} described a reduction of the $\mathrm{2D}$-Toda hierarchy called Toeplitz lattice. This hierarchy turns out to be equivalent to the one originally described by Ablowitz and Ladik \cite{AL} using semidiscrete zero-curvature equations. In this paper we obtain the original semidiscrete zero-curvature equations starting directly from the Toeplitz lattice and we generalize these computations to the matrix case. This generalization lead us to the semidiscrete zero-curvature equations for the non-abelian (or multicomponent) version of the Ablowitz-Ladik equations \cite{GI}. In this way we extend the link between biorthogonal polynomials on the unit circle and the Ablowitz-Ladik hierarchy to the matrix case. \\
MSC classes: 05E35, 37K10.
\end{abstract}

\section{Introduction.}
The Ablowitz-Ladik hierarchy has been introduced in 1975 \cite{AL} as a spatial discretization of the $\mathrm{AKNS}$ hierarchy. As described by Suris in \cite{Su} Ablowitz and Ladik replaced the celebrated Zakharov-Shabat spectral problem 
\begin{gather*}
\begin{cases}
	\partial_x \Psi=L\Psi\quad\\\\
	\partial_\tau \Psi=M\Psi
\end{cases}
\end{gather*}
with a discretized version of it; namely
\begin{gather*}
\begin{cases}
	\Psi_{k+1}=L_k \Psi_k\\\\
	\partial_\tau\Psi_k=M_k\Psi_k.
\end{cases}
\end{gather*}
Here $\Psi$ and $\Psi_k$ are two-component vectors while $L,M,L_k$ and $M_k$ are $2\times2$ matrices. In particular
\begin{equation*}
	L:=\begin{pmatrix}
					z & x\\
					y & -z
			\end{pmatrix}
\end{equation*}
and
\begin{equation*}
	L_k:=\begin{pmatrix}
				z & x_k\\
				y_k & z^{-1}
			\end{pmatrix}.
\end{equation*}
We can consider, as usual, the periodic case ($k\in\bZ_n$), the infinite case ($k\in\bZ$) or the semi-infinite case ($k\in\bN$). The standard zero-curvature equations for the $\mathrm{AKNS}$ hierarchy are replaced with the semidiscrete zero-curvature equations
\begin{equation*}\label{semidiscrete}
	\partial_\tau L_k=M_{k+1}L_k-L_kM_k.
\end{equation*} 
As an example, one of the most important equation of this hierarchy is the discrete complexified version of the nonlinear Schr\"odinger equations
\begin{equation*}
	\begin{cases}
		\partial_\tau x_k=x_{k+1}-2x_k+x_{k-1}-x_ky_k(x_{k+1}+x_{k-1})\\\\
		\partial_\tau y_k=-y_{k+1}+2y_k-y_{k-1}+x_ky_k(y_{k+1}+y_{k-1}).
	\end{cases}
\end{equation*}

Quite recently different authors (see \cite{AVM} and \cite{Ne}) worked on the link between biorthogonal polynomials on the unit circle and the semi-infinite Ablowitz-Ladik hierarchy. This is an analogue of the celebrated link between Toda hierarchy and orthogonal polynomials on the real line. Many articles have been written about orthogonal polynomials and the Toda hierarchy and it is almost impossible to recall all of them. Let us just mention the seminal paper of Moser \cite{Mo} and the monography \cite{Te}; there interested readers can find many references. The more recent articles \cite{HM}, \cite{AVM1} and \cite{AVM} should be also cited. In these articles the connection between the theory of orthogonal polynomials and the Sato's theory of infinite Grassmannians appeared for the first time.
In particular in \cite{AVM} the case of Toda and Ablowitz-Ladik hierarchies are treated in a similar way as reductions of $\mathrm{2D}$-Toda. In this way it's given a clear and unified explanation of the role played by orthogonal polynomials in the theory of integrable equations. As noted by the authors their approach is quite different from the original one; actually in their paper \cite{AVM} they always speak about Toeplitz lattice and the coincidence with Ablowitz-Ladik hierarchy is just stated in the introduction. Nevertheless, as explained in the second section of this paper, it is very easy to deduce semidiscrete zero-curvature equations starting from Adler-van Moerbeke's equations. Having this approach in mind we addressed a question arising from the following facts:
\begin{itemize}
	\item Time evolution for orthogonal polynomials on the real line leads to the Toda hierarchy.
	\item Time evolution for biorthogonal polynomials on the unit circle leads to the Ablowitz-Ladik hierarchy.
	\item Time evolution for matrix orthogonal polynomials on the real line leads to the non-abelian Toda hierarchy (see for instance \cite{BG} and references therein; our approach is close to \cite{M2}).
\end{itemize}
\emph{What about time evolution for matrix biorthogonal polynomials on the unit circle?} 
\\

In other words our goal was to replace the question mark in the table below with the corresponding hierarchy.\\
\begin{center}
{
\offinterlineskip
\tabskip=0pt
\halign{ 
\vrule height2.75ex depth1.25ex width 0.6pt #\tabskip=1em &
\hfil #\hfil &\vrule # & #\,\hfil &\vrule # &
\hfil #\hfil &#\vrule width 0.6pt \tabskip=0pt\cr
\noalign{\hrule height 0.6pt}
& \omit &&\omit Orthogonal polynomials on $\bR$ && \omit Biorthogonal polynomials on $S^1$ & \cr
\noalign{\hrule}
& scalar case && \qquad\qquad\quad Toda && Ablowitz-Ladik &\cr
& matrix case && \qquad \:non-Abelian Toda && ? &\cr
\noalign{\hrule height 0.6pt}
}}
\end{center}

In this article we prove that the relevant hierarchy is the non-Abelian version of the Ablowitz-Ladik hierarchy. This hierarchy has been already studied by different authors since 1983 (see \cite{GI}, \cite{AOT} and \cite{TUW}) but, at our best knowledge, the connection with matrix biorthogonal polynomials on the unit circle was never established before. We also remark that this hierarchy is usually called matrix, vector or multicomponent Ablowitz-Ladik. Here we prefer to call it non-Abelian to stress the analogy with Toda. In our setting, instead of having one Lax operator $L_k$ of size $2\times 2$, we have two Lax operators $\cal L_k^l$ and $\cal L_k^r$ of size $2n\times 2n$. These operators depend on matrices $x_k^l,x_k^r,y_k^l$ and $y_k^r$ of size $n\times n$. Semidiscrete zero curvature equations are given by
\begin{gather*}
		\partial_{\tau}\cal L_k^l=\cal M_{k+1}^l\cal L_{k}^l-\cal L_n^l\cal M_{k}^l\\
		\partial_{\tau}\cal L_k^r=\cal L_{k}^r\cal M_{k+1}^r-\cal M_{k}^r\cal L_k^r.
	\end{gather*}
	where also $\cal M_{k}^l$ and $\cal M_{k}^r$ are block matrices.
For instance we have, in this hierarchy, two versions of the non-Abelian complexified discrete nonlinear Schr\"odinger:
\begin{equation*}
	\begin{cases}
		\partial_\tau x_k^l=x_{k+1}^l-2x_k^l+x_{k-1}^l-x_{k+1}^ly_k^rx_k^l-x_k^ly_k^rx_{k-1}^l\\\\
		\partial_\tau y_k^r=-y_{k+1}^r+2y_k^r-y_{k-1}^r+y_{k+1}^rx_k^ly_k^r+y_k^rx_k^ly_{k-1}^r
	\end{cases}
\end{equation*}
\begin{equation*}
	\begin{cases}
		\partial_\tau x_k^r=x_{k+1}^r-2x_k^r+x_{k-1}^r-x_{k-1}^ry_k^lx_k^r-x_k^ry_k^lx_{k+1}^r\\\\
		\partial_\tau y_k^l=-y_{k+1}^l+2y_k^l-y_{k-1}^l+y_{k-1}^lx_k^ry_k^l+y_k^lx_k^ry_{k+1}^l
	\end{cases}
\end{equation*}
(see \cite{APT} for a review about these equations).\\
Sections are organized as follows:
\begin{itemize}
	\item In the second section we recall some basic facts about 2D-Toda and the connection with biorthogonal polynomials. We use the approach developed in \cite{AVM}.
	\item Third section starts with the description of the Toeplitz lattice (see \cite{AVM}) and shows how to deduce semidiscrete zero-curvature equations for the Ablowitz-Ladik hierarchy.
	\item In the fourth section we extend the Toeplitz lattice to the case of block Toeplitz matrices.
	\item Fifth section gives recursion relations for matrix biorthogonal polynomials on the unit circle; these formulas slightly generalize formulas contained in \cite{M} and \cite {As} for matrix orthogonal polynomials on the unit circle.
	\item In the sixth section we derive block semidiscrete zero-curvature equations defining the non-Abelian Ablowitz-Ladik hierarchy. As an example we write the non-abelian analogue of the discrete nonlinear Schr\"odinger. 
\end{itemize}

\section{2D-Toda; linearization and biorthogonal polynomials.}

In this section we recall some basic facts about the 2D-Toda hierarchy as presented in \cite{T}. Moreover we describe the connection with biorthogonal polynomials as originally presented in \cite{AVM}.
We are interested in the semi-infinite case; we start denoting with $\Lambda$ the shift matrix
$$\Lambda:=(\delta_{i+1,j})_{i,j\geq 0}.$$
For the transpose we use the notation $\Lambda^{\mathrm T}=\Lambda^{-1}.$
Then we define two Lax matrices
\begin{gather*}
\begin{cases}
	L_1:=\Lambda+\sum_{i\leq 0}a_i^{(1)}\Lambda^i\\\\
	L_2:=a_{-1}^{(2)}\Lambda^{-1}+\sum_{i\geq 0}a_i^{(2)}\Lambda^{i}
\end{cases}
\end{gather*}
where $\lbrace a_i^{(s)},s=1,2\rbrace$ are some diagonal matrices.
2D-Toda equations, expressed in Lax form, arise as compatibility conditions for the following Zakarov-Shabat spectral problem:
\begin{gather*}
\begin{cases}
	L_1\Psi_1=z\Psi_1\\\\
	L_2^T\Psi^*_2=z^{-1}\Psi^*_2\\\\
	\partial_{t_n}\Psi_1=(L_1^n)_+\Psi_1\\\\
	\partial_{t_n}\Psi_2^*=-(L_1^n)_+^{\mathrm T}\Psi_2^*\\\\
	\partial_{s_n}\Psi_1=(L_2^n)_-\Psi_1\\\\
	\partial_{s_n}\Psi_2^*=-(L_2^n)_-^{\mathrm T}\Psi_2^*.
\end{cases}
\end{gather*}
Here we introduced two infinite sets of times $\lbrace t_i,i\geq 0\rbrace$ and $\lbrace s_i,i\geq 0\rbrace$. We denoted with $N_+$ the upper triangular part of a matrix $N$ (including the main diagonal) and with $N_-$ the lower triangular part (excluding the main diagonal).
$\Psi_1$ and $\Psi_2^*$ are semi-infinte column vectors of type
\begin{gather*}
		\Psi_1(z)=(\Psi_{1,0}(z),\Psi_{1,1}(z),\ldots)^{\mathrm T}\\
 		\Psi_2^*(z)=(\Psi_{2,0}^*(z),\Psi_{2,1}^*(z),\ldots)^{\mathrm T}.
\end{gather*}
For every $k$ the two expressions $e^{-\xi(t,z)}\Psi_{1,k}(z)$ and $e^{-\xi(s,z)}\Psi_{2,k}^*(z^{-1})$ are polynomials in $z$ of order $k$.
Lax equations are written as 
$$\partial_{t_n}L_i=\big[(L_1^n)_+,L_i\big]\quad \partial_{s_n}L_i=\big[(L_2^n)_-,L_i\big],\quad i=1,2.$$
2D-Toda equations can be linearized as explained in \cite{T}. We start with an initial value matrix $M(0,0)=\lbrace M_{ij}(0,0)\rbrace_{i,j\geq 0}$ and we define its time evolution through the equation
\begin{equation*}
	M(t;s):=\exp\Big(\xi(t,\Lambda)\Big)M(0,0)\exp\Big(-\xi(s,\Lambda^{-1})\Big).
\end{equation*}
We assume that there exist a factorization
\begin{gather*}
	M(0,0)=S_1(0,0)^{-1}S_2(0,0).
\end{gather*}
Here $S_1$ is lower triangular while $S_2$ is upper triangular. We assume that both $S_1$ and $S_2$ have non zero elements on the main diagonal and we normalize them in such a way that every element on the main diagonal of $S_1$ is equal to $1$.
Moreover we consider values of $t$ and $s$ for which we can write
\begin{gather}\label{LU}
	M(t,s)=S_1(t,s)^{-1}S_2(t,s)
\end{gather}
with $S_1$ and $S_2$ having the same properties as above. It can be proven that such factorization exists if and only if all the principal minors of $M(t,s)$ do not vanish. In particular this condition is satisfied when $M(t,s)$ is the matrix of the moments of a positive-definite measure. Now we denote with $\chi(z)$ the infinite vector $\chi(z):=(1,z,z^2,\ldots)^{\mathrm T}.$
Wave vectors for $\mathrm{2D}$-Toda and Lax matrices are constructed in the following way.
\begin{thm}[\cite{T}]
	The wave vectors
	\begin{gather*} 
		\Psi_1(z):=\exp(\xi(t,z))S_1\chi(z)\\
		\Psi_2^*(z):=\exp(-\xi(s,z^{-1}))(S_2^{-1})^{\mathrm T}\chi(z^{-1}).
	\end{gather*}
	and the two Lax operators $L_1:=S_1\Lambda S_1^{-1}$ and $L_2:=S_2\Lambda^{-1} S_2^{-1}$ satisfy the $\mathrm{2D}$-Toda Zakharov-Shabat spectral problem.
\end{thm}
\prf
We just sketch the proof and make reference to the article \cite{T}. It is clear that the matrix $M(t,s)$ satisfies differential equations
\begin{gather*}
	\partial_{t_i}M=\Lambda^iM\\
	\partial_{s_i}M=-M\Lambda^{-i}.
\end{gather*}
Then it is easy to deduce Sato's equations
\begin{gather*}
		\partial_{t_n}S_1=-(L_1^n)_-S_1 \\
		\partial_{t_n}S_2=(L_1^n)_+S_2	\\
		\partial_{s_n}S_1=(L_2^n)_-S_1	\\
		\partial_{s_n}S_2=-(L_2^n)_+S_2.
\end{gather*}
and Zakharov-Shabat's equations can be deduced from the expression of wave vectors in terms of $S_1$ and $S_2$. 
\endprf

The last thing we need is the link between the factorization of $M$ and biorthogonal polynomials. We introduce a bilinear pairing on the space of polynomials in $z$ defining
\begin{equation*}
	<z^i,z^j>_M:=M_{ij}.
\end{equation*}
The following proposition is a direct consequence of (\ref{LU}).
\begin{prop}[\cite{AVM}]
	\begin{gather*}
	q^{(1)}=(q^{(1)}_i)_{i\geq 0}:=S_1\chi(z)\\
	q^{(2)}=(q^{(2)}_i)_{i\geq 0}:=(S_2^{-1})^T\chi(z)
\end{gather*}
are biorthonormal polynomials with respect to the pairing $<,>_M$; i.e.
\begin{equation*}
	<q^{(1)}_i,q^{(2)}_j>_M=\delta_{ij}\quad\forall\:i,j\in\bN.
\end{equation*}
\end{prop} 

\section{From the Toeplitz lattice hierarchy to the semidiscrete zero-curvature equations for the Ablowitz-Ladik hierarchy.}
In this section we briefly recall the reduction from $\mathrm{2D}$-Toda to the Toeplitz lattice as described in \cite{AVM}. Then we show how the Ablowitz-Ladik equations are easily obtained from the Toeplitz lattice.
Suppose that our initial value $M(0,0)$ is a Toeplitz matrix; i.e. we have
\begin{equation*}
	M(0,0)=T(\gamma)=\begin{pmatrix}
											\gamma^{(0)} & \gamma^{(-1)} & \gamma^{(-2)} & \ldots\\
											&&&\\
											\gamma^{(1)} & \gamma^{(0)} & \gamma^{(-1)} & \ldots\\
											&&&\\
											\gamma^{(2)} & \gamma^{(1)} & \gamma^{(0)} & \ldots\\
											&&&\\
											\vdots & \vdots & \vdots & \ddots
										\end{pmatrix}
\end{equation*}
for some formal power series $\gamma(z)=\sum_{n\in\bZ}\gamma^{(n)}z^n$. Since $\Lambda=T(z^{-1})$ is an upper triangular Toeplitz matrix it follows easily (see for instance \cite{BoSi}) that
\begin{gather*}
	M(t,s)=\exp\Big(\xi\big(t,\Lambda\big)\Big)M(0,0)\exp\Big(-\xi\big(s,\Lambda^{-1}\big)\Big)=\\
	T\Big(\exp\big(\xi(t,z^{-1})\big)\gamma(z)\exp\big(-\xi(s,z)\big)\Big).
\end{gather*}
This means that Toeplitz form in conserved along $\mathrm{2D}$-Toda flow, hence we are dealing with a reduction of it. This reduction is called Toeplitz lattice in \cite{AVM}. In that article the authors noticed, in the introduction, that this is nothing but the Ablowitz-Ladik hierarchy. Now we will describe how to obtain the original formulation of the Ablowitz-Ladik equations starting from Adler-van Moerbeke's formulation.
The key observation is that, in this case, the bilinear pairing $<p,q>_M$ between two arbitrary polynomials is given by 
\begin{equation*}
	<p,q>_M=\oint p(z)\gamma(z)q^*(z)\ddz.
\end{equation*}
Here the symbol of integration means that we are taking the residue of the formal series $p(z)\gamma(z)q^*(z)$ and $q^*(z)=q(z^{-1})$. In other words $q^{(1)}_i$ and $q^{(2)}_j$ are nothing but orthonormal polynomials on the unit circle. 
We also define monic biorthogonal polynomials
\begin{gather*}
	p^{(1)}=(p^{(1)}_i)_{i\geq 0}:=S_1\chi(z)\\
	p^{(2)}=(p^{(2)}_i)_{i\geq 0}:=h(S_2^{-1})^T\chi(z)
\end{gather*}
with $h=\mathrm{diag}(h_0,h_1,h_2,\ldots)$ some diagonal matrix.
Given an arbitrary polynomial $q(z)$ of degree $n$ we define its reversed polynomial $\tilde q(z):=z^nq^*(z)$ and reflection coefficients
\begin{gather*}
	x_{n}:=p^{(1)}_n(0)\quad y_n:=p^{(2)}_n(0).
\end{gather*} 
We can state the standard recursion relation associated to biorthogonal polynomials on the unit circle. The equation below was already known to Szeg\"o \cite{Sz}. It has been used in \cite{H} in relation with the theory of integrable equations.
\begin{prop}\label{scalarrecursion}
	The following recursion relation holds:
	\begin{equation}\label{recursion}
		\begin{pmatrix}
			p_{n+1}^{(1)}(z)\\
			\tilde p_{n+1}^{(2)}(z)
		\end{pmatrix}=\cal L_n\begin{pmatrix}
			p_{n}^{(1)}(z)\\
			\tilde p_{n}^{(2)}(z)
			\end{pmatrix}=\begin{pmatrix}
												z 			& x_{n+1}\\
												zy_{n+1}& 1
											\end{pmatrix}\begin{pmatrix}
			p_{n}^{(1)}(z)\\
			\tilde p_{n}^{(2)}(z)
		\end{pmatrix}
	\end{equation}
\end{prop}
Using this recursion relation Adler and van Moerbeke in \cite{AVM} wrote the peculiar form of Lax operators for the Toeplitz reduction.
\begin{prop}[\cite{AVM}]
	Lax operators of Toeplitz lattice are of the following form:
	\begin{gather*}
		h^{-1}L_1h=\begin{pmatrix} -x_1y_0 & 1-x_1y_1 & 0 & \ldots & \ldots\\
														 -x_2y_0 & -x_2y_1 & 1-x_2y_2 & 0 & \ldots\\
														 -x_3y_0 & -x_3y_1 & -x_3y_2 & 1-x_3y_3& 0\\
														 \vdots & \vdots & \vdots & \vdots & \vdots\\
														 \vdots & \vdots & \vdots & \vdots & \vdots
							\end{pmatrix}\\\\
						L_2=\begin{pmatrix}-x_0y_1 & -x_0y_2 & -x_0y_3 & \ldots & \ldots\\
														 1-x_1y_1 & -x_1y_2 & -x_1y_3 & \ldots & \ldots\\
														 0 & 1-x_2y_2 & -x_2y_3 & \ldots& \ldots\\
														 \vdots & 0 & 1-x_3y_3 & \ldots & \ldots\\
														 \vdots & \vdots & 0 & \ldots & \ldots
								\end{pmatrix}.
	\end{gather*}
\end{prop}
From this proposition it follows easily the following corollary.
\begin{corol}
	Reflection coefficients satisfy the following equation
	$$\frac{h_{n+1}}{h_n}=1-x_{n+1}y_{n+1}.$$
\end{corol}

Now we can state the theorem relating the Toeplitz lattice to the original form of the Ablowitz-Ladik hierarchy.

\begin{thm}\label{MEvolution}
	The Toeplitz lattice flow can be written in the form
	\begin{gather}
		\partial_{t_i}\begin{pmatrix}
				p_{n}^{(1)}(z)\\
				\tilde p_{n}^{(2)}(z)
				\end{pmatrix}=\cal M_{t_i,n}\begin{pmatrix}
				p_{n}^{(1)}(z)\\
				\tilde p_{n}^{(2)}(z)
			\end{pmatrix}\label{timeevolution1}\\
		\partial_{s_i}\begin{pmatrix}
				p_{n}^{(1)}(z)\\
				\tilde p_{n}^{(2)}(z)
				\end{pmatrix}=\cal M_{s_i,n}\begin{pmatrix}
				p_{n}^{(1)}(z)\\
				\tilde p_{n}^{(2)}(z)
			\end{pmatrix}\label{timeevolution2}
\end{gather}
for some matrices $\cal M_{t_i,n},\cal M_{s_i,n}$ depending on $\lbrace x_j,y_j,z\rbrace$.
\end{thm}
\prf We prove it for the set of times denoted with $t$.\\
We denote $d([z])=\mathrm{diag}(1,z,z^2,z^3,\ldots)$. We have the identities
\begin{gather*}
	\Psi_1(z)=\exp(\xi(t,z))p^{(1)}(z)\label{wave1}\\
	\Psi_2^*(z)=h^{-1}d([z^{-1}])\exp(-\xi(s,z^{-1}))\tilde p^{(2)}(z)\label{wave2}
\end{gather*}
giving the following time evolution for orthogonal polynomials
\begin{gather*}
		\partial_{t_i}p^{(1)}(z)=-z^ip^{(1)}(z)+(L_1^i)_+p^{(1)}(z)\\
		\partial_{t_i}\tilde p^{(2)}(z)=
		-hd([z])(L_1^i)_{++}^Th^{-1}d([z^{-1}])\tilde p^{(2)}(z).
\end{gather*}
Here $(L_1^i)_{++}$ denotes the strictly upper diagonal part of $L_1^n$. The formulas above are obtained from a straightforward computation and using the fact, proven in \cite{AVM}, that 
$$\partial_{t_i}\log(h_n)=(L_1^i)_{nn}.$$
Hence we have that, for every $k$, $\partial_{t_i}p^{(1)}_k$ is a linear combination of $\lbrace p^{(1)}_k,p^{(1)}_{k+1},p^{(1)}_{k+2},\ldots\rbrace$ with coefficients in $\bC[x_j,y_j]$.
In the same way, for every $k$, $\partial_{t_i}\tilde p^{(2)}_k$ is a linear combination of $\lbrace\tilde p^{(2)}_k,\tilde p^{(2)}_{k-1},\tilde p^{(2)}_{k-2},\ldots\rbrace$ with coefficients in $\bC[x_j,y_j]$.
Using the recursion relation (\ref{recursion}) and its inverse
\begin{equation*}\label{recursioninverse}
		\begin{pmatrix}
			p_{n}^{(1)}(z)\\
			\tilde p_{n}^{(2)}(z)
		\end{pmatrix}=\cal L_n^{-1}\begin{pmatrix}
			p_{n+1}^{(1)}(z)\\
			\tilde p_{n+1}^{(2)}(z)
			\end{pmatrix}=\frac{h_n}{h_{n+1}}
								\begin{pmatrix}
									z^{-1} & -z^{-1}x_{n+1}\\
									-y_{n+1} & 1		
								\end{pmatrix}\begin{pmatrix}
			p_{n+1}^{(1)}(z)\\
			\tilde p_{n+1}^{(2)}(z)
		\end{pmatrix}
\end{equation*}
we can obtain the desired matrices $\cal M_{t_i,n}.$
\endprf\\
\begin{corol}[Ablowitz-Ladik semidiscrete zero-curvature equations]
 The matrices $\cal L_n$ satisfy the following time evolution
 \begin{gather}
		\partial_{t_i}\cal L_n=\cal M_{t_i,n+1}\cal L_{n}-\cal L_n\cal M_{t_i,n}\\
		\partial_{s_i}\cal L_n=\cal M_{s_i,n+1}\cal L_{n}-\cal L_n\cal M_{s_i,n}.
	\end{gather}
\end{corol}
\prf
These equations are nothing but the compatibility conditions of recursion relation (\ref{recursion}) with time evolution (\ref{timeevolution1}) and (\ref{timeevolution2}).
\endprf\\
\begin{rem}
	Actually our Lax operator $\cal L_n$ is slightly different from the Lax operator $L_n$ used in \cite{AL} by Ablowitz and Ladik and written in the introduction above. Nevertheless, as shown in \cite{MEKL}, these two Lax operators are linked through a simple change of spectral parameter.
\end{rem}

\begin{ex}[The first flows; discrete nonlinear Schr\"odinger.]
	The first matrices $\cal M_{t_i,n}$ and $\cal M_{s_i,n}$ are easily computed. We have
	\begin{gather*}
	\partial_{t_1}p^{(1)}_k=-zp^{(1)}_k-x_{k+1}y_kp^{(1)}_k+p^{(1)}_{k+1}=-x_{k+1}y_kp^{(1)}_k+x_{k+1}\tilde p^{(2)}_k\\
	\partial_{t_1}\tilde p^{(2)}_k=-z\frac{h_{k+1}}{h_k}\tilde p^{(2)}_{k-1}=zy_kp^{(1)}_k-z\tilde p^{(2)}_k
\end{gather*}
giving immediately
\begin{equation*}
	\cal M_{t_1,k}=\begin{pmatrix}
										-x_{k+1}y_k & x_{k+1}\\
										z y_k       & -z
									\end{pmatrix}
\end{equation*} 
	An analogue computation for $s_1$ can be easily done obtaining
\begin{equation*}
		\cal M_{s_1,k}=\begin{pmatrix}
										z^{-1} & -z^{-1}x_k\\
										-y_{k+1}       & x_ky_{k+1}
									\end{pmatrix}.
\end{equation*}
We can already write, with $t_1$ and $s_1$, the well known integrable discretization of nonlinear Schr\"odinger.
We just need two more trivial rescaling times introduced with substitutions 
\begin{gather*}
		p^{(1)}\longmapsto \exp(t_0)p^{(1)}\\
		\tilde p^{(2)}\longmapsto \exp(-s_0)\tilde p^{(2)}
\end{gather*}
and corresponding to matrices
\begin{gather*}
	\cal M_{t_0,k}:=\begin{pmatrix} 1 & 0\\ 0 & 0\end{pmatrix}\\
	\cal M_{s_0,k}:=\begin{pmatrix} 0 & 0\\ 0 & -1\end{pmatrix}.
\end{gather*}
Now we can construct the matrix
\begin{gather*}
		\cal M_{\tau,k}=\cal M_{t_1,k}+\cal M_{s_1,k}-\cal M_{t_0,k}-\cal M_{s_0,k}=\\\\
		\begin{pmatrix}
				z^{-1}-1-x_{k+1}y_k & x_{k+1}-z^{-1}x_k\\
				zy_k-y_{k+1}        & x_ky_{k+1}+1-z
		\end{pmatrix}
\end{gather*}
associated to the time $\tau=t_1+s_1-t_0-s_0$ so that semidiscrete zero-curvature equation
$$\partial_\tau\cal L_k=\cal M_{\tau,k+1}\cal L_{k}-\cal L_k\cal M_{\tau,k}$$
is equivalent to the system
\begin{equation}\label{dNLS}
	\begin{cases}
		\partial_\tau x_k=x_{k+1}-2x_k+x_{k-1}-x_ky_k(x_{k+1}+x_{k-1})\\\\
		\partial_\tau y_k=-y_{k+1}+2y_k-y_{k-1}+x_ky_k(y_{k+1}+y_{k-1}).
	\end{cases}
\end{equation}
This is exactly the complexified version of the discrete nonlinear Schr\"odinger. Rescaling $\tau\mapsto i\tau$ and imposing $y_k=\pm x^*_k$ we obtain
\begin{equation}\label{dNLS2}
	-i\partial_\tau x_k=x_{k+1}-2x_k+x_{k-1}\mp \|x_k\|^2(x_{k+1}+x_{k-1}).
\end{equation}
\end{ex}

\section{Toda flow for block Toeplitz matrices and the related Lax operators.}

Now we generalize the Toeplitz lattice's equations to the block case. We start with a matrix-valued formal series
$$\gamma(z)=\sum_{k\in\bZ}\gamma^{(k)}z^k.$$
Here every element $\gamma^{(k)}$ is a $n\times n$ matrix. Then we define its time evolution as
$$\gamma(t,s;z):=\exp\big(-\xi(s,z^{-1}\Id)\big)\gamma(z)\exp\big(\xi(t,z\Id)\big).$$
where $\Id$ is the $n\times n$ identity matrix.
Differently from the scalar case we don't consider just one Toeplitz matrix but the two block Toeplitz matrices, right and left, given by
\begin{gather*}
	T^r(\gamma):=\begin{pmatrix}
											\gamma^{(0)} & \gamma^{(-1)} & \gamma^{(-2)} & \ldots\\
											&&&\\
											\gamma^{(1)} & \gamma^{(0)} & \gamma^{(-1)} & \ldots\\
											&&&\\
											\gamma^{(2)} & \gamma^{(1)} & \gamma^{(0)} & \ldots\\
											&&&\\
											\vdots & \vdots & \vdots & \ddots
										\end{pmatrix}\\\\
	T^l(\gamma):=\begin{pmatrix}
											\gamma^{(0)} & \gamma^{(1)} & \gamma^{(2)} & \ldots\\
											&&&\\
											\gamma^{(-1)} & \gamma^{(0)} & \gamma^{(1)} & \ldots\\
											&&&\\
											\gamma^{(-2)} & \gamma^{(-1)} & \gamma^{(0)} & \ldots\\
											&&&\\
											\vdots & \vdots & \vdots & \ddots
										\end{pmatrix}
\end{gather*}
In this way we obtain the following linear time evolution for our block Toeplitz matrices (in the following we will omit the symbol $\gamma$):
\begin{gather}
	\partial_{t_i}T^l=\Lambda^iT^l\quad \partial_{t_i}T^r=T^r\Lambda^{-i}\\
	\partial_{s_i}T^l=-T^l\Lambda^{-i}\quad \partial_{s_i}T^r=-\Lambda^iT^r
\end{gather}
where, in this case, we have $\Lambda=T^r(z^{-1}\Id)$.
Then we assume that there exist two factorizations
\begin{gather*}
	T^l=S_1^{-1}S_2 \quad T^r=Z_2Z_1^{-1}.
\end{gather*}
Here $S_1,Z_2$ are block-lower triangular while $S_2,Z_1$ are block-upper triangular. We assume that all these matrices have non degenerate blocks on the main diagonal (i.e. these blocks must have non zero determinants). Normalizations are chosen in such a way that every element on the main  block-diagonal of $S_1$ and $Z_2$ is equal to the identity matrix $\Id$. As we did before we assume that these conditions hold when every time is equal to $0$ and the we consider just the values of $t$ and $s$ for which these conditions still hold.
In the matrix case we can define two bilinear pairings given by the following definition.
\begin{defin}
	$$<P,Q>_r:=\oint P^*(z)\gamma(z)Q(z)\ddz\quad <P,Q>_l:=\oint P(z)\gamma(z)Q^*(z)\ddz$$
	where $P$ and $Q$ are two arbitrary matrix polynomials and $P^*(z):=(P(z^{-1}))^{\mathrm T}$
\end{defin}
Our two factorizations give exactly biorthonormal polynomials for $<,>_r$ and $<,>_l$. In the following we denote $\chi(z):=(\Id ,z\Id,z^2\Id,z^3\Id,\ldots)^{\mathrm T}.$
\begin{prop}\label{dressingbiorthonormal}
	\begin{gather}
		Q^{(1)l}:=\begin{pmatrix} Q^{(1)l}_0\\Q^{(1)l}_1\\ \vdots\end{pmatrix}=S_1\chi(z)\\
		Q^{(2)l}:=\begin{pmatrix} Q^{(2)l}_0\\Q^{(2)l}_1\\ \vdots\end{pmatrix}=(S_2^{-1})^{\mathrm T}\chi(z)\\
		Q^{(1)r}:=\begin{pmatrix} Q^{(1)r}_0 & Q^{(1)r}_1 & \ldots\end{pmatrix}=\chi(z)^{\mathrm T} Z_1\\
		Q^{(2)r}:=\begin{pmatrix} Q^{(2)r}_0 & Q^{(2)r}_1 & \ldots\end{pmatrix}=\chi(z)^{\mathrm T} (Z_2^{-1})^{\mathrm T}
	\end{gather}
	are the biorthonormal polynomials associated to the pairing $<,>_l$ and $<,>_r$. Hence for every $i,j$ we have
	$$<Q^{(1)l}_i,Q^{(2)l}_j>_l=\delta_{ij}\quad <Q^{(2)r}_i,Q^{(1)r}_j>_r=\delta_{ij}$$
\end{prop}
\prf
 	We just prove, as an example, the proposition for the right polynomials. On the other hand the one for left polynomials is identical to the usual proof for $\mathrm{2D}$-Toda. We have
 	\begin{gather*}
 		\Bigg(<Q^{(2)r}_i,Q^{(1)r}_j>_r\Bigg)_{i,j\geq 0}=\Bigg(\sum_{k,l\geq 0}(Z_2^{-1})_{ki}<z^k\Id,z^l\Id>_r(Z_1)_{lj}\Bigg)=\\\\
 		Z_2^{-1}T^lZ_1=\Id\Longleftrightarrow T^l=Z_2Z^{-1}.
 	\end{gather*}
 	(it should be noted that, in this case, subscripts of type $(Z_1)_{ij}$ denote the block in position $(i,j)$ and not the element $(i,j)$.)
\endprf\\
Now we can write the related Sato's equations for $S_i$ and $Z_i$. It is convenient to introduce the following Lax operators.
\begin{defin}
	\begin{gather}
		L_1:=S_1\Lambda S_1^{-1} \quad L_2:=S_2\Lambda^{-1}S_2^{-1}\\
		R_1:=Z_1^{-1}\Lambda^{-1}Z_1 \quad R_2:=Z_2^{-1}\Lambda Z_2.
	\end{gather}
\end{defin}
 	
\begin{prop}
	The following Sato's equations are satisfied.
	\begin{gather}
		\partial_{t_n}S_1=-(L_1^n)_-S_1 \quad \partial_{t_n}Z_1=-Z_1(R_1^n)_+\\
		\partial_{t_n}S_2=(L_1^n)_+S_2	\quad \partial_{t_n}Z_2=Z_2(R_1^n)_-\\
		\partial_{s_n}S_1=(L_2^n)_-S_1	\quad \partial_{s_n}Z_1=Z_1(R_2^n)_+\\
		\partial_{s_n}S_2=-(L_2^n)_+S_2	\quad \partial_{s_n}Z_2=-Z_2(R_2^n)_-
	\end{gather}
\end{prop} 	
 	\prf
 		We will just prove, as an example, the equations involving $t$-derivative of $Z_1$ and $Z_2$.\\
 		We assume, as an ansatz, that we have 
 		\begin{gather*}
 			\partial_{t_n}Z_1=Z_1A\\
 			\partial_{t_n}Z_2=Z_2B.
 		\end{gather*}
 		for some matrices $A$ and $B$.
 		Then exploiting time evolution of $T^r$ we can write
 		\begin{gather*}
 			T^r\Lambda^{-n}=\partial_{t_n}T^r=\partial_{t_n}(Z_2Z_1^{-1})=\\
 			Z_2BZ_1^{-1}-Z_2Z_1^{-1}Z_1AZ_1^{-1}=Z_2(B-A)Z_1^{-1}
 		\end{gather*}
 		hence we must have $(B-A)Z_1^{-1}=Z_1^{-1}\Lambda^{-n}$. We rewrite it as 
 		$$B-A=Z_1^{-1}\Lambda^{-n}Z_1=R_1^n.$$ Moreover we have that $B$ must be strictly block-lower triangular and $A$ block-upper triangular. This is because $Z_1$ is block-upper triangular and $Z_2$ is block-lower triangular with constant entries on the main diagonal. Hence we conclude that $A=-(R_1^n)_+$ and $B=(R_1^n)_-$.
 		\endprf\\
Now it's just a matter of trivial computations to write down the corresponding Lax equations for $L_i$ and $R_i$.

\begin{prop}
	The following Lax equations are satisfied:
	\begin{gather}
		\partial_{t_n}L_i=\Big[(L_1^n)_+,L_i\Big] \quad \partial_{t_n}R_i=\Big[R_i,(R_1^n)_-\Big]\label{Lax1}\\
		\partial_{s_n}L_i=\Big[(L_2^n)_-,L_i\Big] \quad \partial_{s_n}R_i=\Big[R_i,(R_2^n)_+\Big].\label{Lax2}
	\end{gather}
\end{prop}
The definition of our Lax operators will give us eigenvalue equations for suitably defined wave vectors.
\begin{defin}
	\begin{gather}
		\Psi_1(z):=\exp(\xi(t,z\Id))S_1\chi(z) \\	\Phi_1(z):=\exp(\xi(t,z\Id))\Big[\chi(z)\Big]^{\mathrm T}Z_1\\
		\Psi_2^*(z):=\exp(-\xi(s,z^{-1}\Id))(S_2^{-1})^{\mathrm T}\chi(z^{-1}) \\ 	\Phi_2^*(z):=\exp(-\xi(s,z^{-1}\Id))\chi(z^{-1})^{\mathrm T}(Z_2^{-1})^{\mathrm T}.
	\end{gather}
\end{defin}
\begin{prop}
	The following equations hold true:
	\begin{gather}
		L_1\Psi_1(z)=z\Psi_1(z) \quad \Phi_1(z)R_1=z\Phi_1(z)\label{Eigenvalue1}\\
		L_2^{\mathrm T}\Psi_2^*(z)=z^{-1}\Psi_2^*(z) \quad \Phi_2^*(z)R_2^{\mathrm T}=z^{-1}\Phi_2^*(z).\label{Eigenvalue2}
	\end{gather} 
\end{prop}
\prf We will just prove the last equation, all the other ones are proved in a similar way.
	From the very definition we have
	\begin{gather*}
		\Phi_2^*(z)R_2^{\mathrm T}=z^{-1}\Id\Phi_2^*(z)
		\Longleftrightarrow [\chi(z^{-1})]^{\mathrm T}(Z_2^{-1})^{\mathrm T}R_2^{\mathrm T}=
		z^{-1}[\chi(z^{-1})]^{\mathrm T}(Z_2^{-1})^{\mathrm T}=\\
		[\chi(z^{-1})]^T\Lambda^{-1}(Z_2^{-1})^{\mathrm T} 
		\Longleftrightarrow R_2^{\mathrm T}=Z_2^{\mathrm T}\Lambda^{-1}(Z_2^{-1})^{\mathrm T}
		\Longleftrightarrow R_2=Z_2^{-1}\Lambda Z_2
	\end{gather*}
\endprf\\
The proof of the following proposition is straightforward.
\begin{prop}
	The Lax equations (\ref{Lax1}) and (\ref{Lax2}) are the compatibility conditions of the eigenvalue equations (\ref{Eigenvalue1}) and (\ref{Eigenvalue2}) with the following equations:
	\begin{gather}
		\partial_{t_n}\Psi_1=(L_1^n)_+\Psi_1 \quad \partial_{t_n}\Phi_1=\Phi_1(R_1^n)_-\label{ZS1}\\
		\partial_{s_n}\Psi_1=(L_2^n)_-\Psi_1 \quad \partial_{s_n}\Phi_1=\Phi_1(R_2^n)_+\label{ZS2}\\
		\partial_{t_n}\Psi_2^*=-(L_{1+}^n)^{\mathrm T}\Psi_2^* \quad 
		\partial_{t_n}\Phi_2^*=-\Phi_2^*(R_{1-}^n)^{\mathrm T}\label{ZS3}\\
		\partial_{s_n}\Psi_2^*=-(L_{2-}^n)^{\mathrm T}\Psi_2^* \quad 
		\partial_{s_n}\Phi_2^*=-\Phi_2^*(R_{2+}^n)^{\mathrm T}\label{ZS4}.
	\end{gather}
\end{prop}
\begin{rem}
Actually our Lax equations as well as equations (\ref{Eigenvalue1})-(\ref{ZS4}) can be deduced from the equations of multicomponent 2D-Toda \cite{UT}.
\end{rem}

\section{Recursion relations for matrix biorthogonal polynomials on the unit circle.}

In order to generalize the scalar theory we have to construct an analogue of the recursion relation given by Proposition \ref{scalarrecursion}. Recursion relations for matrix orthogonal polynomial on the unit circle are already known, see \cite{M} and \cite{As}. Here we slightly generalize to the case of matrix biorthogonal polynomials on the unit circle.\\
We define the following important $n\times n$ matrices:
\begin{defin}
	$$h_N^r:=\SC(T_{N+1}^r)\quad h_N^l:=\SC(T_{N+1}^l)$$
\end{defin}
where $\SC$ denote the $n\times n$ Schur complement of a block matrix with respect to the upper left block (see for instance \cite{FZ}). For example we have
\begin{equation*}
	\SC(T_{N+1}^r)=
	\gamma^{(0)}-\begin{pmatrix}\gamma^{(N)}&\ldots&\ldots&\gamma^{(1)}\end{pmatrix}T_N^{-r}
	\begin{pmatrix}\gamma^{(-N)}\\\ldots\\\ldots\\\gamma^{(-1)}\end{pmatrix}
\end{equation*}
(here and below $T_N^{-r}:=(T_N^{r})^{-1}$ and similarly for $h_N^r,h_N^l$ and $T_N^l$).
\begin{prop}
	Monic biorthogonal polynomials such that
	$$<P^{(2)r}_k,P^{(1)r}_j>_r=\delta_{kj}h_k^r\quad\quad<P^{(1)l}_k,P^{(2)l}_j>_l=\delta_{kj}h_k^l.$$
	are given by the following formulas:
	$$P^{(1)r}_N=\SC\begin{pmatrix} \gamma^{(0)} & \ldots & \ldots & \gamma^{(-N+1)} & \gamma^{(-N)}\\
													\ldots & \ldots & \ldots & \ldots & \ldots\\
													\ldots & \ldots & \ldots & \ldots & \ldots\\
													\gamma^{(N-1)} & \ldots & \ldots & \gamma^{(0)} & \gamma^{(-1)}\\
													\Id            & z\Id     & \ldots & z^{N-1}\Id & z^N\Id
						\end{pmatrix}$$
	$$(P^{(2)r}_N)^T=\SC\begin{pmatrix} \gamma^{(0)} & \ldots & \ldots & \gamma^{(-N+1)} & \Id\\
													\ldots & \ldots & \ldots & \ldots & \ldots\\
													\ldots & \ldots & \ldots & \ldots & \ldots\\
													\gamma^{(N-1)} & \ldots & \ldots & \gamma^{(0)} & z^{N-1}\Id\\
													\gamma^{(N)}   & \ldots & \ldots & \gamma^{(1)} & z^N\Id
						\end{pmatrix}$$
	$$P^{(1)l_N}=\SC\begin{pmatrix} \gamma^{(0)} & \ldots & \ldots & \gamma^{(N-1)} & \Id\\
													\ldots & \ldots & \ldots & \ldots & \ldots\\
													\ldots & \ldots & \ldots & \ldots & \ldots\\
													\gamma^{(-N+1)} & \ldots & \ldots & \gamma^{(0)} & z^{N-1}\Id\\
													\gamma^{(-N)}   & \ldots & \ldots & \gamma^{(-1)} & z^N\Id
						\end{pmatrix}$$
	$$(P^{(2)l}_N)^T=\SC\begin{pmatrix} \gamma^{(0)} & \ldots & \ldots & \gamma^{(N-1)} & \gamma^{(1)}\\
													\ldots & \ldots & \ldots & \ldots & \ldots\\
													\ldots & \ldots & \ldots & \ldots & \ldots\\
													\gamma^{(-N+1)} & \ldots & \ldots & \gamma^{(0)} & \gamma^{(N-1)}\\
													\Id            & z\Id     & \ldots & z^{N-1}\Id & z^N\Id
						\end{pmatrix}.$$
\end{prop}
\prf
We will just prove the first formula, the second one is proved similarly.\\
First of all we have $\forall \:0\leq m\leq N-1$
\begin{gather*}
	<z^m\Id,P_N^{(1)r}(z)>_r=
	\oint z^{-m}\gamma(z)
	\Bigg(z^N\Id-\begin{pmatrix}\Id&\ldots&\ldots&z^{N-1}\Id\end{pmatrix}T_N^{-r}
	\begin{pmatrix}\gamma^{(-N)}\\ \ldots \\ \ldots \\ \gamma^{(-1)}\end{pmatrix}\Bigg)\ddz\\
	=\gamma^{(m-N)}-\gamma^{(m-N)}=0.
\end{gather*}
In the same way $\forall\: 0\leq m\leq N-1$
\begin{gather*}
	<P_N^{(2)r}(z),z^m\Id>_r=
	\oint \Bigg(z^{-N}\Id-\begin{pmatrix}\gamma^{(N)}&\ldots&\ldots&\gamma^{(1)}\end{pmatrix}T_N^{-r}
	\begin{pmatrix}\Id\\ \ldots \\ \ldots \\ z^{-N+1}\Id\end{pmatrix}\Bigg)\gamma(z)z^m\ddz\\
	=\gamma^{(N-m)}-\gamma^{(N-m)}=0.
\end{gather*}
Finally
\begin{gather*}
	<P_N^{(2)r}(z),P_N^{(1)r}(z)>_r=<z^N\Id,P^{(1)r}_N>_r=\\
	\gamma^{(0)}-
	\begin{pmatrix}\gamma^{(N)}&\ldots&\ldots&\gamma^{(1)}\end{pmatrix}T_N^{-r}
	\begin{pmatrix}\gamma^{(-N)}\\ \ldots \\ \ldots \\ \gamma^{(-1)}\end{pmatrix}=h^r_N.
\end{gather*}
This completes the proof of the first formula, the second one is proved similarly.\endprf\\
\begin{rem}\label{normalization}
	Note that, imposing 
	\begin{gather*}
		Q^{(1)l}_k:=P^{(1)l}_k \quad Q^{(2)l}_k:=(h_k^{-l})^TP^{(2)l}_k\\
		Q^{(1)r}_k:=P^{(1)r}_k(h_k^{-r}) \quad Q^{(2)r}_k:=P^{(2)r}_k.
	\end{gather*}
	we obtain biorthonormal polynomials.
\end{rem}
Now we will write a long list of relations among this polynomials and reflection coefficients.
Given any matrix polynomial $Q(z)$ of degree $n$ we define the associated reversed polynomial as
$$\tilde Q(z)=z^nQ^*(z).$$
The reflection coefficients are defined as follows:
\begin{gather*}
	x_{N}^l:=P^{(1)l}_N(0) \quad x_N^r:=P^{(1)r}_N(0)\\
	y_{N}^l:=(P^{(2)l}_N(0))^{\mathrm T} \quad y_{N}^r:=(P^{(2)r}_N(0))^{\mathrm T}.
\end{gather*}
\begin{prop}
	The following formulas hold true:
	\begin{gather}
		P^{(1)l}_{N+1}-zP_N^{(1)l}=x_{N+1}^l\tilde P_N^{(2)r}\label{r1}\\
		\tilde P^{(2)r}_{N+1}-\tilde P_N^{(2)r}=zy_{N+1}^rP_N^{(1)l}\label{r2}\\
		P^{(1)r}_{N+1}-zP_N^{(1)r}=\tilde P_N^{(2)l}x_{N+1}^r\label{r3}\\
		\tilde P^{(2)l}_{N+1}-\tilde P_N^{(2)l}=zP_N^{(1)r}y_{N+1}^l\label{r4}\\
		P^{(1)r}_{N+1}=zP^{(1)r}_N(\Id-y^l_{N+1}x^r_{N+1})+\tilde P_{N+1}^{(2)l}x^r_{N+1}\label{r7}\\
		P^{(1)l}_{N+1}=z(\Id-x^l_{N+1}y^r_{N+1})P^{(1)l}_N+x^l_{N+1}\tilde P_{N+1}^{(2)l}\label{r8}\\
		\tilde P^{(2)r}_{N+1}=(\Id-y^r_{N+1}x^l_{N+1})\tilde P^{(2)r}_{N}+y^r_{N+1}P^{(1)l}_{N+1}\label{r9}\\
		\tilde P^{(2)l}_{N+1}=\tilde P^{(2)l}_{N}(\Id-x^r_{N+1}y^l_{N+1})+P^{(1)r}_{N+1}y^l_{N+1}\label{r10}\\
		x^l_Nh^r_N=h^l_Nx^r_N\label{r5}\\
		y^r_Nh^l_N=h^r_Ny^l_N\label{r6}\\
		h_N^{-r}h_{N+1}^r=\Id-y^l_{N+1}x^r_{N+1}\label{r11}\\
		h_{N+1}^{l}h_{N}^{-l}=\Id-x^l_{N+1}y^r_{N+1}.\label{r12}
	\end{gather}
\end{prop}
\prf
The first four formulas are proved observing, for instance for the first case, that $\forall\: 1\leq i\leq N$ we have
$$0=<P^{(1)l}_{N+1}-zP^{(1)l}_N,z^i \Id>_l=<\tilde P^{(2)r}_N,z^i\Id>_l$$
so that $P^{(1)l}_{N+1}-zP^{(1)l}_N$ and $\tilde P^{(2)r}_N$ must be proportional. Setting $z=0$ you also find the constant of proportionality.
In particular, when proving (\ref{r2}) and (\ref{r4}), we find a formula and then we have to take the reversed one.
(\ref{r7}) is proved substituting (\ref{r4}) into (\ref{r3}) and similarly for (\ref{r8}),(\ref{r9}),(\ref{r10}).\\
(\ref{r5}) and (\ref{r6}) are proven respectively observing that we have
$$<\tilde P^{(1)l}_N,P^{(1)r}_N>_r=<P^{(1)l}_N,\tilde P^{(1)r}_N>_l$$
and
$$<P^{(2)r}_N,\tilde P^{(2)l}_N>_r=<\tilde P^{(2)r}_N,P^{(2)l}_N>_l$$
and then doing explicit computations.
Finally (\ref{r11}) is obtained rewriting (\ref{r7}) as
$$\dfrac{P^{(1)r}_{N+1}}{z^{N+1}}=\dfrac{P^{(1)r}_N}{z^N}(\Id-y^l_{N+1}x^r_{N+1})+(P^{(2)l}_{N+1})^*x^r_{N+1},$$
multiplying from the left for $P^{(1)l}_N\gamma$ and then taking the residue.
(\ref{r12}) is proved similarly.
\endprf\\
Now we define two sets of block matrices $\lbrace\cal L_N^r\rbrace_{N\geq 0}$ and $\lbrace\cal L_N^l\rbrace_{N\geq 0}$. They will have, in the matrix case, the same role played by $\lbrace\cal L_n\rbrace_{n\geq 0}$ in the scalar case.
\begin{defin}
	\begin{gather}
		\cal L_N^l:=
		\begin{pmatrix}
										z\Id & x^l_{N+1}\\
										zy^r_{N+1} & \Id
		\end{pmatrix}\label{sLax1}\\
		\cal L_N^r:=
		\begin{pmatrix}
			z\Id & zy^l_{N+1}\\
			x^r_{N+1}& \Id
		\end{pmatrix}.\label{sLax2}
	\end{gather}
\end{defin}
\begin{corol}
	The following block matrices recursion relations are satisfied
	\begin{gather}
		\begin{pmatrix}
			P^{(1)l}_{N+1}\\
			\tilde P^{(2)r}_{N+1}
		\end{pmatrix}=\cal L_N^l\begin{pmatrix}
			P^{(1)l}_{N}\\
			\tilde P^{(2)r}_{N}
		\end{pmatrix}
		\label{blockrecursion1}\\
		\begin{pmatrix}
			P^{(1)r}_{N+1}&
			\tilde P^{(2)l}_{N+1}
		\end{pmatrix}=\begin{pmatrix}
			P^{(1)r}_{N}&
			\tilde P^{(2)l}_{N}
		\end{pmatrix}\cal L_N^r.
		\label{blockrecursion2}
	\end{gather}
\end{corol}
\prf These are nothing but (\ref{r1}),(\ref{r2}),(\ref{r3}) and (\ref{r4}).\endprf\\

\section{Explicit expressions for Lax operators and related semidiscrete zero-curvature equations.}

Using our recursion relations we want to find explicit expressions for $L_i$ and $R_i$ in terms of our reflection coefficients
$x_k^l,x_k^r,y_k^l,y_k^r$. First of all we underline a remarkable symmetry that will allow us to reduce the amount of our computations. Doing the following three substitutions
\begin{gather*}
	z\mapsto z^{-1}\\
	t\mapsto -s\\
	s\mapsto -t
\end{gather*}
we obtain immediately the following proposition.
\begin{prop}
	Under the symmetry above the dressings, the orthogonal polynomials, the Lax operators and the reflection coefficients change as follows:
	\begin{gather*}
		T^r\mapsto T^l \quad T^l\mapsto T^r\\
		S_1\mapsto Z_2^{-1} \quad S_2\mapsto Z_1^{-1}\\
		L_1\mapsto R_2 \quad L_2\mapsto R_1\\
		Q^{(1)l}\mapsto (Q^{(2)r})^* \quad Q^{(2)l}\mapsto (Q^{(1)r})^*\\
		P^{(1)l}\mapsto (P^{(2)r})^* \quad P^{(2)l}\mapsto (P^{(2)r})^*\\
		x^l\mapsto y^r \quad y^l\mapsto x^r\\
		h^l_k\mapsto h^r_k.
	\end{gather*}
\end{prop}
Hence we can write just the left theory and we will have the right one as well. Actually every computation made above for the right theory can be deduced from the left theory and this symmetry which will be called in the sequel $t-s$ symmetry.
In the theorem below the symbol $\prod_{j=N+2}^{M-}$ means that the terms in the product must be written from the smallest index to the biggest, going from left to right. The symbol $\prod_{j=N+2}^{M+}$ means that the product must be taken in the opposite direction.
\begin{thm}[Lax operators for block Toeplitz lattice]
	Lax operators $L_i$ and $R_i$ are expressed in terms of reflection coefficients according to the following formulas:\\ $\forall\:N>M\geq-1$
	\begin{gather}
		(L_1)_{N,M+1}=-x^l_{N+1}\Big(\prod_{j=N+2}^{M-}(\Id-y^r_jx^l_j)\Big)y^r_{M+1}\label{h1}\\
		(R_2)_{N,M+1}=-y^r_{N+1}\Big(\prod_{j=N+2}^{M-}(\Id-x^l_jy^r_j)\Big)x^l_{M+1}\label{h2}\\
		(L_2)_{M+1,N}=-h^{-l}_{M+1}x^r_{M+1}\Big(\prod_{j=N+2}^{M+}(\Id-y^l_jx^r_j)\Big)y^l_{N+1}h^{l}_{N}\label{h3}\\
		(R_1)_{M+1,N}=-h^{-r}_{M+1}y^l_{M+1}\Big(\prod_{j=N+2}^{M+}(\Id-x^r_jy^l_j)\Big)x^r_{N+1}h^{r}_{N}.\label{h4}
	\end{gather}
	Moreover 
	\begin{gather}
		(L_1)_{N,N+1}=(R_2)_{N,N+1}=\Id\label{h5}\\
		(L_2)_{N+1,N}=h^l_{N+1}h^{-l}_{N}\quad
		(R_1)_{N+1,N}=h^r_{N+1}h^{-r}_{N}.\label{h6}
	\end{gather}
\end{thm}
\prf
	(\ref{h5}) and (\ref{h6}) follow trivially from the expressions of dressings $S_i$, $Z_i$. Infact, because of our normalization, we have that $S_1$ and $Z_2$ are equal to the identity matrix plus a strictly lower triangular matrix which is what stated in (\ref{h5}). (\ref{h6}) is obtained observing that the block on the main diagonal of $S_2$ and $Z_1$ are given respectively by matrices $h_k^l$ and $h_k^r$ as follows from Proposition \ref{dressingbiorthonormal} and Remark \ref{normalization}.
	Now let's begin with (\ref{h1}) and (\ref{h2}); the important point is that we have
$$\Psi_1=\exp(\xi(t,z\Id))P^{(1)l}.$$
Hence as done in \cite{AVM} we can find that $\forall\:N>M\geq-1$ we have
	\begin{equation}\label{h1a}
		(L_1)_{N,M+1}=-x^l_{N+1}h^r_Nh^{-r}_{M+1}y^r_{M+1}.
	\end{equation}
Infact $\forall\:N>M\geq -1$
	\begin{gather*}
		<P^{(1)l}_{N+1}-zP^{(1)l}_N,P^{(2)l}_{M+1}-zP^{(2)l}_M>_l=-<zP^{(1)l}_N,P^{(2)l}_{M+1}>_l=\\
		-<P^{(1)l}_{N+1}+\ldots+(L_1)_{N,M+1}P^{(1)l}_{M+1}+\ldots,P^{(2)l}_{M+1}>_l=-(L_1)_{N,M+1}h^l_{M+1}.
	\end{gather*}
	On the other hand, using recursion relations, I also have $\forall N\geq M\geq -1$
	\begin{gather*}
		<P^{(1)l}_{N+1}-zP^{(1)l}_N,P^{(2)l}_{M+1}-zP^{(2)l}_M>_l=<x^l_{N+1}\tilde P^{(2)r}_N,(y^l_{M+1})^T\tilde P^{(1)r}_M>_l=\\
		x^l_{N+1}\Bigg(\oint z^{N-M}(P^{(2)r}_N)^*\gamma(z)P^{(1)r}_M\ddz\Bigg) y^l_{M+1}=\\
		x^l_{N+1}<P^{(2)r}_N,z^{N-M}P^{(1)r}_M>_r y^l_{M+1}=x^l_{N+1}h^r_N y^l_{M+1}
	\end{gather*}
and comparing them we find (\ref{h1a}).
Now we use $t-s$ symmetry to simplify this expression. We obtain
$$(R_2)_{N,M+1}=-y^r_{N+1}h^l_Nh^{-l}_{M+1}x^l_{M+1}.$$
Now we apply several times the recursion (\ref{r12}) to the piece $h^l_Nh^{-l}_{M+1}$ and we get (\ref{h2}). (\ref{h1}) is obtained using $t-s$ symmetry.
For (\ref{h3}) and (\ref{h4}) we start defining $\tilde R_1$ such that $z P^{(1)r}=P^{(1)r}\tilde R_1$; then we will have
$R_1=h^{r}\tilde R_1h^{-r}$ and computations for $\tilde R_1$ is carried on similarly as for $L_1$.
\endprf\\
\begin{rem}
	Our equations (\ref{h1}),(\ref{h2}),(\ref{h3}) and (\ref{h4}) extend to the matrix biorthogonal setting the equations written in \cite{As} for matrix orthogonal polynomials (see equation (4.2),(4.3)). In that article properties of $M$ are applied to study some problems in computational mathematics (multivariate time series analysis and multichannel signal processing) and no relation is established with the Lax theory and integrable systems.
\end{rem}
The theorem above describes the block-analogue of Toeplitz lattice; we are now in the position to prove the analogue of Theorem \ref{MEvolution}.

\begin{thm}
	Block Toeplitz lattice flow can be written in the form
	\begin{gather}
		\partial_{t_i/s_i}\begin{pmatrix}
				P_{N}^{(1)l}(z)\\
				\tilde P_{N}^{(2)r}(z)
				\end{pmatrix}=\cal M_{t_i/s_i,N}^l\begin{pmatrix}
				P_{N}^{(1)l}(z)\\
				\tilde P_{N}^{(2)r}(z)
				\end{pmatrix}\label{blockt1}\\
		\partial_{t_i/s_i}\begin{pmatrix}
				P_{N}^{(1)r}(z) & \tilde P_{N}^{(2)l}(z)
				\end{pmatrix}=\begin{pmatrix}
				P_{N}^{(1)r}(z) & \tilde P_{N}^{(2)l}(z)
				\end{pmatrix}\cal M_{t_i/s_i,N}^r\label{blockt2}
\end{gather}
for some block matrices $\cal M_{t_i,N}^r,\cal M_{s_i,N}^r,\cal M_{t_i,N}^l,\cal M_{s_i,N}^l$ depending on the matrices $\lbrace x_j^l,y_j^l,x_j^r,y_j^r\rbrace$ and the spectral parameter $z$.
\end{thm}
\prf
As we did for the scalar case we prove it just for $t$ times. The relevant equations linking biorthogonal polynomials with wave vectors are
\begin{gather*}
	\Psi_1(z)=\exp(\xi(t,z\Id))P^{(1)l}(z)\\
	\Phi^*_2(z)=\exp(-\xi(s,z^{-1}\Id))(\tilde P^{(2)r})^{\mathrm T} d([z^{-1}])\\
	\Phi_1(z)=\exp(\xi(t,z\Id))P^{(1)r}h^{-r}\\
	\Psi_2^*(z)=\exp(-\xi(s,z^{-1}\Id))(h^{-l})^{\mathrm T}d([z^{-1}])(\tilde P^{(2)l})^{\mathrm T}.
\end{gather*}
Trivial computations give the following time evolution:
\begin{gather*}
	\partial_{t_n}P^{(1)l}=(L_1^n)_+P^{(1)l}-z^{n} P^{(1)l}\\
	\partial_{t_n}\tilde P^{(2)r}=-d([z])(R_1^n)_-d([z^{-1}])\tilde P^{(2)r}\\
	\partial_{t_n}P^{(1)r}=P^{(1)r}(h^{-r}(R_1^n)_-h^r+h^{-r}(\partial_{t_n}h^r)-z^n\Id)\\
	\partial_{t_n}\tilde P^{(2)l}=\tilde P^{(2)l}(-h^{-l}(L_1^n)_+h^l+h^{-l}(\partial_{t_n}h^l)).
\end{gather*}
the last two can be simplified giving
\begin{gather}
	\partial_{t_n}P^{(1)l}=(L_1^n)_+P^{(1)l}-z^{n}\Id P^{(1)l}\label{t1}\\
	\partial_{t_n}\tilde P^{(2)r}=-d([z])(R_1^n)_-d([z^{-1}])\tilde P^{(2)r}\label{t2}\\
	\partial_{t_n}P^{(1)r}=P^{(1)r}(h^{-r}(R_1^n)_{--}h^r-z^n\Id)\label{t3}\\
	\partial_{t_n}\tilde P^{(2)l}=-\tilde P^{(2)l}\Big(d([z^{-1}])h^{-l}(L_1^n)_{++}h^ld([z])\Big).\label{t4}
\end{gather}
where $(R_1^n)_{--}$ means the lower triangular part including the main diagonal and $(L_1^n)_{++}$ means the strictly upper diagonal part. This simplification can be obtained evaluating the terms $h^{-r}(\partial_{t_n}h^r)$ and $h^{-l}(\partial_{t_n}h^l))$ through Sato's equations. Also they can be obtained observing that $P^{(1)r}_N$ and $P^{(2)l}_N$ are monic so that the derivative of the leading term is equal to $0$.
Then the proof is obtained as we did in the scalar case using forward and backward recursion relations (\ref{r1}),(\ref{r3}),(\ref{r9}) and (\ref{r10}).
\endprf
\begin{corol}[Non-abelian AL semidiscrete zero-curvature equations]
 The matrices $\cal L_n^r$ and $\cal L_n^l$ satisfy the following time evolution
 \begin{gather}
		\partial_{t_i/s_i}\cal L_n^l=\cal M_{t_i/s_i,n+1}^l\cal L_{n}^l-\cal L_n^l\cal M_{t_i/s_i,n}^l\\
		\partial_{t_i/s_i}\cal L_n^r=\cal L_{n}^r\cal M_{t_i/s_i,n+1}^r-\cal M_{t_i/s_i,n}^r\cal L_n^r.
	\end{gather}
\end{corol}
\prf
These equations are nothing but compatibility conditions of recursion relations (\ref{blockrecursion1}) and (\ref{blockrecursion2}) with time evolution (\ref{blockt1}) and (\ref{blockt2}).
\endprf\\

\begin{rem}
It should be noticed that, with respect to the equations originally written in \cite{GI}, here we have two coupled non-abelian Ablowitz-Ladik equations.
\end{rem}

\begin{ex}[The first flows; non-abelian analogue of discrete nonlinear Schr\"odinger]
	As we did for the scalar case we will compute the first matrices $\cal M_{t_1/s_1,k}^{r/l}$ and use them to construct the non-abelian version of discrete nonlinear Schr\"odinger. We start with $\cal M_{t_1,k}^l$; (\ref{t1}) gives us immediately
\begin{gather*}
	\partial_{t_1}P^{(1)l}_k=P^{(1)l}_{k+1}-x^l_{k+1}y^r_kP^{(1)l}_k-zP^{(1)l}_k=\\
	zP^{(1)l}_k+x^l_{k+1}\tilde P^{(2)r}_k-x^l_{k+1}y^r_{k}P^{(1)l}_k-zP^{(1)l}_k=-x^l_{k+1}y^r_{k}P^{(1)l}_k+x^l_{k+1}\tilde P^{(2)r}_k
\end{gather*}
while we obtain immediately from (\ref{t2}) that
$$\partial_{t_1}\tilde P^{(2)r}_k=-z h^r_kh^{-r}_{k-1}\tilde P^{(2)r}_{k-1}.$$
Then we use recursion relation (\ref{r9}) combined with
$$h^r_kh^{-r}_{k-1}=(\Id-y^r_kx^l_k)$$
(this one comes from recursion relation (\ref{r5}) combined with $t-s$ symmetry)
to arrive to
$$\partial_{t_1}\tilde P^{(2)r}_k=zy_k^{r}P^{(1)l}_k-zP^{(2)r}_k.$$
These computations give us
\begin{equation}
	\cal M_{t_1,k}^l=\begin{pmatrix}
										-x^l_{k+1}y^r_{k} & x^l_{k+1}\\\\
										zy_k^r & -z\Id
									\end{pmatrix}.
\end{equation}
Exploiting $t-s$ symmetry we can write immediately 
\begin{equation}
	\cal M_{s_1,k}^l=\begin{pmatrix}
										z^{-1}\Id & -z^{-1}x^l_{k}\\\\
										-y_{k+1}^r & y^{r}_{k+1}x^l_{k}
									\end{pmatrix}.
\end{equation}
Analogue computations for $\cal M_{t_1}^r$ gives us

\begin{gather*}
	\partial_{t_1}P^{(1)r}_k=P^{(1)r}_{k+1}-P^{(1)r}_ky^l_kx^r_{k+1}-zP^{1(r)}_k=\\
	zP^{1(r)}_k+\tilde P^{(2)l}_kx^r_{k+1}-P^{(1)r}_ky^l_kx^r_{k+1}-zP^{1(r)}_k=\tilde P^{(2)l}_kx^r_{k+1}-P^{(1)r}_ky^l_kx^r_{k+1}
\end{gather*}
and
\begin{gather*}
	\partial_{t_1}\tilde P^{(2)l}_k=-z\tilde P^{(2)l}_{k-1}(h^{-l}_{k-1}h^{l}_{k})=
	-\tilde P^{(2)l}_kz+P^{(1)r}_kzy^l_k
\end{gather*}
(here we started from (\ref{t3}) and (\ref{t4}) and we used recursion relations (\ref{r7}),(\ref{r10}),(\ref{r6}), the last one combined with $t-s$ symmetry).
Then we arrive to
\begin{equation}
	\cal M_{t_1,k}^r=\begin{pmatrix}
										-y^l_{k}x^r_{k+1} & zy_k^l\\\\
										x^r_{k+1} & -z\Id
									\end{pmatrix}
\end{equation}
and using again $t-s$ symmetry we also get
\begin{equation}
	\cal M_{s_1,k}^r=\begin{pmatrix}
										z^{-1}\Id & -y^l_{k+1}\\\\
										-z^{-1}x_{k}^r & x^r_{k}y^l_{k+1}
									\end{pmatrix}.
\end{equation}
As we did for the scalar case we introduce times $t_0$ and $s_0$ that give matrices
\begin{gather}
	\cal M_{t_0,k}^{r/l}=\begin{pmatrix}
													\Id & 0\\
													0   & 0
												\end{pmatrix}\\
	\cal M_{s_0,k}^{r/l}=\begin{pmatrix}
													0 & 0\\
													0   & -\Id
												\end{pmatrix}.							
\end{gather}
Then we construct the matrices
\begin{gather*}
		\cal M_{\tau,k}^l=\cal M_{t_1,k}^l+\cal M_{s_1,k}^l-\cal M_{t_0,k}^l-\cal M_{s_0,k}^l=\\\\
		\begin{pmatrix}
				z^{-1}\Id-\Id-x_{k+1}^ly_k^r & x_{k+1}^l-z^{-1}x_k^l\\\\
				zy_k^r-y_{k+1}^r        & y_{k+1}^rx_k^l+\Id-z\Id
		\end{pmatrix}\\\\
		\cal M_{\tau,k}^r=\cal M_{t_1,k}^r+\cal M_{s_1,k}^r-\cal M_{t_0,k}^r-\cal M_{s_0,k}^r=\\\\
		\begin{pmatrix}
				z^{-1}\Id-\Id-y_{k}^lx_{k+1}^r & zy^l_k-y^l_{k+1}\\\\
				x^r_{k+1}-z^{-1}x^r_k        & x_{k}^ry_{k+1}^l-z\Id+\Id
		\end{pmatrix}
\end{gather*}
associate to the time $\tau=t_1+s_1-t_0-s_0$. Semidiscrete zero-curvature equations
\begin{gather*}
	\partial_\tau\cal L_k^l=\cal M_{\tau,k+1}^l\cal L_{k}^l-\cal L_k^l\cal M_{\tau,k}^l\\
	\partial_\tau\cal L_k^r=\cal L_{k}^r\cal M_{\tau,k+1}^r-\cal M_{\tau,k}^r\cal L_k^r
\end{gather*}
are equivalent to the systems
\begin{equation}\label{mNLSa}
	\begin{cases}
		\partial_\tau x_k^l=x_{k+1}^l-2x_k^l+x_{k-1}^l-x_{k+1}^ly_k^rx_k^l-x_k^ly_k^rx_{k-1}^l\\\\
		\partial_\tau y_k^r=-y_{k+1}^r+2y_k^r-y_{k-1}^r+y_{k+1}^rx_k^ly_k^r+y_k^rx_k^ly_{k-1}^r
	\end{cases}
\end{equation}
\begin{equation}\label{mNLSb}
	\begin{cases}
		\partial_\tau x_k^r=x_{k+1}^r-2x_k^r+x_{k-1}^r-x_{k-1}^ry_k^lx_k^r-x_k^ry_k^lx_{k+1}^r\\\\
		\partial_\tau y_k^l=-y_{k+1}^l+2y_k^l-y_{k-1}^l+y_{k-1}^lx_k^ry_k^l+y_k^lx_k^ry_{k+1}^l.
	\end{cases}
\end{equation}
Note that both of them are equivalent to the discrete matrix NLS as written, for instance, in \cite{APT}. Using (\ref{mNLSa}) and (\ref{mNLSb}) together we perform the reduction to the hermitian case in a different way from \cite{APT}. First of all we rescale $\tau\mapsto i\tau$ and then we impose 
\begin{gather*}
	y_k^r=\pm(x_k^r)^*\\
	y_k^l=\pm(x_k^l)^*
\end{gather*}
	Note that this reduction (with the sign plus) corresponds to studying the theory of matrix orthogonal polynomials on the unit circle as described in \cite{M} and \cite{As}, hence it is very natural. This reduction gives us the two coupled equations
\begin{equation}
\begin{cases}
	-i\partial_\tau x_k^l=x_{k+1}^l-2x_k^l+x_{k-1}^l\mp x_{k+1}^l(x_k^r)^*x_k^l\mp x_k^l(x_k^r)^*x_{k-1}^l\\\\
	-i\partial_\tau x_k^r=x_{k+1}^r-2x_k^r+x_{k-1}^r\mp x_{k-1}^r(x_k^l)^*x_k^r\mp x_k^r(x_k^l)^*x_{k+1}^r
\end{cases}
\end{equation}
already studied in \cite{AOT} and generalized in \cite{TUW}.
\end{ex}

\begin{rem}
In \cite{MEKL} the authors studied finite gap solutions of the Ablowitz-Ladik hierarchy. It could be interesting to generalize their results to the non-abelian case. We will consider this problem in a subsequent publication.
\end{rem}

\subsection*{Acknowledgments} 

I am very grateful to Professor B. Dubrovin for his constant support and many suggestions given during the preparation of this paper. 
Also I wish to thank Professor Takayuki Tsuchida that gave me useful references to the already existing papers about the non-Abelian Ablowitz-Ladik equations and Professor Manuel Manas for some clarifying discussions about multicomponent 2D-Toda hierarchy and its relationship with this work.
I am also grateful to Professors Pierre van Moerbeke, Mark Adler, Arno Kuijlaars and Walter Van Assche for some general discussions on the theory of orthogonal polynomials and its relations with integrable equations.

\vspace{0.5cm}
 
This work started during my PhD in SISSA and it has been partially supported by the European Science Foundation Programme ``Methods of Integrable Systems, Geometry, Applied Mathematics" (MISGAM), the Marie Curie RTN ``European Network in Geometry, Mathematical Physics and Applications"  (ENIGMA),  and by the Italian Ministry of Universities and Researches (MUR) research grant PRIN 2006 ``Geometric methods in the theory of nonlinear waves and their applications". The support of the Belgian IAP project NOSY (``Nonlinear systems, stochastic processes and statistical mechanics") is gratefully acknowledged.

\end{document}